\documentclass[a4paper,12pt]{amsart}
\usepackage[letterpaper, margin=1.2in]{geometry}
\usepackage{latexsym}
\usepackage{amssymb}
\usepackage{amsmath}
\usepackage{amsfonts}
\usepackage[table]{xcolor}
\usepackage{pgfplots}
\usepackage{color}
\usepackage{multirow}
\usepackage{graphicx}
\usepackage{txfonts,pxfonts} 
\usepackage{tikz}
\usetikzlibrary{calc,arrows}
\usetikzlibrary{calc}
\usepackage{verbatim}
\everymath{\displaystyle}
\newtheorem{thm}{Theorem}[section]

\newtheorem{lem}[thm]{Lemma}

\newtheorem{cor}[thm]{Corollary}

\newtheorem{prop}[thm]{Proposition}
\theoremstyle{definition}

\newtheorem{rem}{Remark}

\def\fph{\mathbb{F}_{\ph}}

\newcommand{\Z}{\mathbb Z}

\newcommand{\Q}{\mathbb Q}

\newcommand{\F}{\mathbb F}

\def\F{\mathbb{F}}

\newcommand{\p}{\mathfrak{p}}
\def\ol{\overline}
\def\al{\alpha}
\def\la{\lambda}

\def\om{\omega}
\def\th{\theta}
\def\md#1{\ \mbox{\rm(mod }{#1})}
\def\nph#1{N_{\ph}(#1)}
\def\npp#1{N_{\ph}^+(#1)}
\def\ph{\phi}

\newcounter{cs}
\stepcounter{cs}
\newcommand{\casos}{\begin{itemize}}
\newcommand{\fcasos}{\end{itemize}\setcounter{cs}{1}}

\newfont{\tit}{cmr12 scaled \magstep3}
\author{ Lhoussain El Fadil }
\email{lhoussain.elfadil@usmba.ac.ma}
\begin{document}
\title[]{On index divisors and  monogenity of certain  number fields defined by trinomials $x^7+ax+b$}
\textcolor[rgb]{1.00,0.00,0.00}{}
\textcolor[rgb]{1.00,0.00,0.00}{}
\begin{abstract}
 
 In this paper  for every number field $K$ generated by a root $\al$ of a trinomial $x^7+ax+b\in \Z[x]$ and for every prime integer $p$, we calculate $\nu_p(i(K))$, the highest power of $p$ dividing the index $i(K)$ of the field $K$. In particular, we calculate the index $i(K)$. When the index of $K$ is not trivial, namely $\nu_p(i(K))\ne q 0$ for some prime integer $p$, then  $K$ is not monogenic. 
\end{abstract}
\maketitle
\vspace{0.3cm}
\section{Introduction}
Let $K$ be a number field of degree $n$, $\Z_K$ its ring of integers, and $d_K$ its absolute discriminant. It is well known that $\Z_K$ is a free abelian group of rank $n$ and by the The fundamental theorem of finite abelian groups, $(\Z_K:\Z[\th])$ for every primitive element $\th\in \Z_K$ of $K$. Let  $ind(\th)=(\Z_K:\Z[\th])$. Then  $ind(\th)$  is called the index of $\th$. The index of the  number field $K$ is defined by $i(K)=GCD((\mathbb{Z}_K:\mathbb{Z}[\theta])\mid K=\mathbb{Q}(\theta) \mbox{ and } \theta\in \Z_K)$. A rational prime integer $p$ dividing $i(K)$ is called a prime common index divisor of $K$. 
The number field $K$ is called { monogenic} if it admits  a $\Z$ basis of type $(1,\th,\ldots,\th^{n-1})$ for  some $\th\in\Z_K$.  Remark that if $\mathbb{Z}_K$ has a power integral basis, then $i(K)=1$. Therefore a field having a prime common index divisor is not monogenic.
Monogenity of number fields is a classical problem of algebraic number theory, going back to Dedekind, Hasse and Hensel, see for instance \cite{Ha,  He} and \cite{G19} for the present state of this area. 
It is called a problem of Hasse to give an arithmetic characterization of those number fields which are monogenic \cite{Ha, He, MNS}.
For any primitive element $\th\in \Z_K$ of $K$,
we denote by 
\[
ind(\th)=(\Z_K:\Z[\th])
\]
the { index of $\th$}, that is 
the index of the $\Z$-module $\Z[\th]$ in the free-$\Z$-module $\Z_K$ of rank $n$.
As it is known \cite{G19}, we have
\[
|\triangle(\th)|=ind(\th)^2\cdot |d_K|
\]
where $\triangle(\th)$ is the discriminant of the minimal polynomial of  $\th$ over $\Q$.

Clearly, $ind(\th)=1$  for some primitive element $\th\in \Z_K$ of $K$ if and only if $(1,\th,\ldots, \th^{n-1})$ is a power integral basis
of $\Z_K$. 
\smallskip

The problem of testing the monogenity of number fields and constructing power integral bases have been intensively studied during the last four decades mainly by Ga\'al, Gy\"ory, Nakahara, Pohst and their collaborators (see for instance \cite{AN, 13a, P}).
 In 1871, Dedekind was the first who gave an example of a number field with non trivial index, he considered the cubic field $K$ generated by a root of $x^3-x^2-2x-8$ and showed that the rational prime $2$ splits completely in $K$ {(\cite[§ 5, page 30]{R})}. According to a well known theorem of Dedekind (\cite[Chapter I, Proposition 8.3]{Neu}), if we suppose that $K$ is monogenic, then we would be able to find a cubic polynomial defining $K$, that splits completely into distinct polynomials of degree $1$ in $\mathbb{F}_2[x]$. Since there is only two distinct polynomials of degree $1$ in $\mathbb{F}_2[x]$, this is impossible. In $1930$, Engstrom  was the first one who  related the prime ideal factorization and the index of a number field of degree less than $8$ \cite{En}. For any number field $K$ of degree $n\leq 7$, he showed that $\nu_p(i(K))$  is explicitly determined by the factorization of $p\mathbb{Z}_K$ into powers of prime ideals of $K$ for every positive rational prime integer  $p\leq n$. This motivated Narkiewicz to ask a very important question, stated as problem 22 in Narkiewicz's book (\cite[Problem 22]{Nar}), which asks for an explicit formula of the highest power $\nu_p(i(K))$ for a given rational prime $p$ dividing $i(K)$.  In \cite{Na}, Nakahara studied the index of non-cyclic but abelian biquadratic number fields. He showed that the field index of such fields is in the set $\{1,2,3,4,6,12\}$. In \cite{GPP}  Ga\'al et al. characterized  the field indices of biquadratic number fields having Galois group $V_4$ and they proved that $i(K)\in\{1,2,3,4,6,12\}$. 
 Recently, many authors are interested on monogenity of number fields defined by trinomials.
In \cite{AK, SK}, Khanduja et al. studied the integral closedness of some number fields defined by trinomials. Their results are refined by Ibarra et al. (see \cite{Smt}).  Recall that the results given in \cite{AK, SK}, can only decide on the integral closedness of $\Z[\al]$, but cannot test whether the field is monogenic or not, where $\al$  is an integral element over $\Z$. In \cite{Lsc, Lac, LPh, LD}, Jones et al. introduced   monogenity of some irreducible polynomials. Recall that  according to {Jones'}  definition, if a polynomial $F(x)$ is monogenic, then $\Q(\al)$ is monogenic, but the converse is not true;  it is possible  that a number field generated by a  root $\al$ of a non monogenic polynomial can be monogenic (see Example \ref{mono}). Therefore {Jones'} and Khanduja's results  cover partially the study of monogenity of number fields defined by trinomials.
Davis and Spearman \cite{DA} studied the index of quartic number fields $K$ generated by a root of 
such a quartic  trinomial $F(x)=x^4+ax+b\in\Z[x]$.
They gave necessary and sufficient conditions on $a$ and $b$ so that a prime $p$ is a common index divisor of $K$ for $p=2,3$. 
Their method is based on the calculation of the $p$-index form of $K$, using $p$-integral bases of $K$.
El Fadil and Ga\'al  \cite{FG4} studied the index of quartic number fields $K$ {generated by a root of a quadratic trinomial of the form  $F(x)=x^4+ax^2+b\in\Z[x]$}.
They gave necessary and sufficient conditions on $a$ and $b$ so that a prime $p$ is a common index divisor of $K$ for every prime integer $p$.
In  \cite{Ga21}, for a sextic number field $K$ defined by a trinomial $F(x)=x^6+ax^3+b\in \Z[x]$, Ga\'al studied the multi-monegenity of $K$; he calculated all possible power integral bases of $K$. In \cite{jnt6}, El Fadil extended Ga\'al's  studies  by providing some cases where $K$ is not monogenic. Also in  \cite{com5}, for every prime integer $p$, El Fadil  gave necessary and sufficient conditions  on $a$ and $b$ so that $p$ is a common index divisor of $K$, where $K$ is a number field defined by  an irreducible trinomial $F(x)=x^5+ax^2+b\in \Z[x]$. 
In this paper,  for a septic number field generated by  a  root of a trinomial $F(x)=x^7 +ax+b\in \Z[x]$, we calculate $\nu_p(i(K))$, the highest power of $p$ dividing the index $i(K)$ of the field $K$. When the index of $K$ is not trivial, namely $\nu_p(i(K))\neq 0$ for some prime integer $p$, then  $K$ is not monogenic. 
	\section{Main Results}
Throughout this section  $K$ is a number field generated by a  root $\al$ of a trinomial $F(x)=x^7+ax+b\in \Z[x]$ and we assume that for every rational prime  integer $p$,  $\nu_p(a)  \le 5$ or  $ \nu_p(b)\le 6$.
  Along this paper, for every integer $a\in \Z$ and a prime integer $p$, let $a_p=\frac{a}{p^{\nu_p(a)}}$.\\
   \smallskip
 
 We start  with the following theorem, which characterizes when is  $\Z[\al]$ integrally closed?
 \begin{thm}\label{intclos}
The ring  $\Z[\al]$ is integrally closed  if and only if the following conditions hold:
\begin{enumerate}
\item
For every prime integer $p$, if $p$ divides $a$ and $b$, then $\nu_p(b)=1$.
\item
 If $2$ divides $b$ and does not  $a$,  then  $a\equiv 3\md4$ and $b \equiv 2\md4$ or $a \equiv 1\md4$ and $b \equiv 0\md4$ .
\item
If $3$ divides $b$ and   $a\equiv -1\md3$, then $\nu_3(b+1+a)=\nu_3(b-1-a)=1$.
\item
If $3$ divides $b$ and   $a\equiv 1\md3$, then   $a\equiv 4\md9$ or $a\equiv 7\md9$ or $\nu_3(b)=1$.
\item
If $7$ divides   $a$ and does not divide $b$, then  $\nu_7(1-a-b^6)=1$.
\item
For every prime integer $p\not\in\{2,3,7\}$, if  $p$ does not divide both $a$ and $b$, then $\nu_p(7^7b^6+6^6a^7)\le 1$.
\end{enumerate}
 \end{thm} 

The following example gives  an infinite family of monogenic  septic number fields defined by  non monogenic trinomials.  In such a way, the results given in \cite{SK, Lsc, Lac, LPh, LD, Smt}, cannot decide on  monogenity of  $K$.\\

{\begin{prop}\label{mono}
		Let  $K$ be the number field generated by a  root $\al$ of $F(x) = x^{7} + 2^u a x + 2^vb  \in \Z[x] $, with $u\ge v-1$, $2\le v\le 6$, GCD$(6,b)=1$, $7$ does not divide $a$, and for every odd prime integer $p$, if $p$  does not divide $b$, then $p^2$ does not divide $7^7b^6+6^6a^7$. Then $F(x)$ is a non monogenic polynomial and $K$ is a monogenic number field.
\end{prop}}
 
   In the remainder of this section, for every prime integer $p$, for every values of  $a$ and $b$ and for every prime integer $p$, we calculate  $\nu_p(i(K))$. If $\nu_p(i(K))\neq 0$ for some prime integer $p$, then $K$ is not monogenic. For every integers $a$ and $b$, let $\triangle=-(6^6a^7+7^7b^6)$ and for every prime integer $p$, let $\triangle_p=\frac{\triangle}{p^{\nu_p(\triangle)}}$.
\begin{thm} \label{npib2}
The following table provides the value of  $\nu_2(i(K))$.
	$$\begin{array}{|c|c|c|}
	\hline
 {\mbox conditions}& \nu_2(i(K))\\
 \hline
	    a\equiv 28\md{32} { \mbox{ and }} b\equiv 0\md{32}&1\\
	    \hline
    a \equiv 112\md{128}{ \mbox{ and }}  b\equiv 0\md{128} &1\\
      \hline
	   a\equiv 1\md8 { \mbox{ and }} b\equiv 2\md{4}&\\
    \nu_2(\triangle) { \mbox{ even and }}\triangle_2\equiv 3\md4&1\\
     \hline
	    a\equiv 3\md8{ \mbox{ and }} b\equiv 4\md8&1\\
	     \hline
	    a\equiv 3\md4{ \mbox{ and }}b\equiv 0\md8&3\\
	          \hline
           (a,b)\in\{(5,2),(5,6),(13,2),(13,14)\}\md{16} &1\\
        \hline
	     {\mbox Otherwise }&0\\
 
   \hline
	     	\end{array}$$
       
		 		In particular, if $\nu_2(i(K))\ge 1$, then  $K$ is not monogenic.
\end{thm}
\begin{thm} \label{npib3}
The following table provides the value of  $\nu_3(i(K))$.
	$$\begin{array}{|c|c|c|}
	\hline
 {\mbox conditions}& \nu_3(i(K))\\
 \hline
	    a\equiv 5\md{9} { \mbox{ and }} b\in\{3,6\} \md{9}&1\\
	    \hline
    a\equiv 8\md{9} { \mbox{ and }} b\equiv 0 \md{9} &2\\
      \hline      
	   a\equiv 2\md9 { \mbox{ and }} b\in\{3,6\}\md{9}&\\
    \nu_3(\triangle)=2k { \mbox{ and }}k\ge 5&1\\
     \hline
	    a\equiv 2\md9 { \mbox{ and }} b\in\{3,6\}\md{9}&\\
    \nu_3(\triangle)=2k+1,\, k\ge 5 { \mbox{ and }}\triangle_3\equiv 1\md3&2\\
        \hline
	     {\mbox Otherwise }&0\\
 
   \hline
	     	\end{array}$$
       
		 		In particular, if $\nu_3(i(K))\ge 1$, then  $K$ is not monogenic.  
 \end{thm}
 \begin{thm}\label{pge5}
     For every prime integer $p\ge 5$ and for every integers $a$ and $b$ such that $F(x)=x^7+ax+b$ is irreducible over $\Q$, $p$ does not divide $i(K)$, where $K$ is the number field defined by $F(x)$.
 \end{thm}
 \begin{cor}
For every integers $a$ and $b$ such that $F(x)=x^7+ax+b$ is irreducible over $\Q$, $i(K)\in\{1, 2, 3, 6, 8, 9, 18, 24, 72\}$, where $K$ is the number field generated by a root of $F(x)$.
\end{cor}
\begin{rem}
\begin{enumerate}
    \item 
The field $K$ can be non monogenic even if the index $i(K)=1$. 
\item
The unique method which allows to test whether $ K $ is monogenic is to calculate the solutions of the index form equation of the field $K$ (see for instance \cite{GR17, G19}).
\end{enumerate}
\end{rem}
\section{A short introduction to prime ideal factorization based on Newton polygons}
\label{intro}

In 1894, Hensel
developed a powerful approach by showing that for every prime integer $p$, the prime ideals of $\Z_K$
lying above   $p$ are in one--one correspondence with
monic irreducible factors of $F(x)$ in $\Q_p[x]$. For every prime ideal
corresponding to any irreducible factor in $\Q_p[x]$, the
ramification index  and the residue degree together are the same as
those of the local field defined  by the associated irreducible factor
\cite{H}. Since then, to factorize $p\Z_K$, we need to factorize $F(x)$ in $\Q_p[x]$.  Newton's polygon techniques can be used to refine the factorization. This is a standard method which is rather technical but
very efficient to apply. We have introduced the corresponding concepts in several former papers.
Here we only give a brief introduction which makes our proofs understandable.
For a detailed description, we refer to Ore's Paper \cite{O} and Guardia, Montes and Nart's paper \cite{GMN}.	
	For every prime integer $p$, let $\nu_p$ be the $p$-adic valuation of $\Q_p$ and $\Z_p$ the ring of $p$-adic integers.  Let $F(x)\in\mathbb{Z}_p[x]$ be a monic polynomial and $\phi\in\mathbb{Z}_p[x]$  a monic lift of an irreducible factor of $\ol{F(x)}$ modulo $p$. Let $F(x)=a_0(x)+a_1(x)\phi(x)+\cdots+a_n(x)\phi(x)^l$ be the $\phi$-expansion of $F(x)$, $\nph{F}$ the $\phi$-Newton polygon of $F(x)$ and $\npp{F}$ its  principal part. Let $\mathbb{F}_{\phi}$ be the field $\mathbb{F}_p[x]/(\overline{\phi})$. For every side $S$ of $\npp{F}$ with length $l$ and initial point $(s,u_s)$,  for every $i=0,\ldots,l$, let  $c_i\in\mathbb{F}_{\phi}$ be the residue coefficient, defined as follows: 
	$$c_{i}=
	\left
	\{\begin{array}{ll} 0,& \mbox{ if } (s+i,{\it u_{s+i}}) \mbox{ lies strictly
			above } S,\\
		\left(\dfrac{a_{s+i}(x)}{p^{{\it u_{s+i}}}}\right)
		\,\,
		\mod{(p,\phi(x))},&\mbox{ if }(s+i,{\it u_{s+i}}) \mbox{ lies on }S.
	\end{array}
	\right.$$
	Let $-\lambda=-h/e$ be the slope of $S$, where $h$ and $e$ are two positive coprime integers. Then  $d=l/e$ is the degree of $S$. Let ${R_{1}(F)(y)}=t_dy^d+t_{d-1}y^{d-1}+\cdots+t_{1}y+t_{0}\in\mathbb{F}_{\phi}[y]$, called  
	the residual polynomial of $F(x)$ associated to the side $S$, where for every $i=0,\dots,d$,  $t_i=c_{ie}$. If ${R_{1}(F)(y)}$ is square free for each side of the polygon $\npp{F}$, then we say that $F(x)$ is $\phi$-regular.\\ 
	Let $\overline{F(x)}=\prod_{i=1}^{r}\ol{\phi_i}^{l_i}$ be the factorization of $F(x)$ into powers of monic irreducible coprime polynomials over $\mathbb{F}_p$, we say that the polynomial $F(x)$ is $p$-{regular} if $F(x)$ is a $\phi_i$-regular polynomial with respect to $p$ for every $i=1,\dots,r$. Let  $N_{\phi_i}^+(F)=S_{i1}+\cdots+S_{ir_i}$ be the $\phi_i$-principal Newton polygon of $F(x)$ with respect to $p$. For every $j=1,\dots,r_i$, let $R_{1_{ij}}(F)(y)=\prod_{s=1}^{s_{ij}}\psi_{ijs}^{a_{ijs}}(y)$ be the factorization of $R_{1_{ij}}(F)(y)$ in $\mathbb{F}_{\phi_i}[y]$, where $R_{1_{ij}}(F)(y)$ is the residual polynomial of $F(x)$ attached to the side $S_{ij}$. Then we have the following  theorem of index of Ore:
	\begin{thm}\label{ore}$($\cite[Theorem 1.7 and Theorem 1.9]{EMN}$)$\\
		Under the above hypothesis, we have the following:
		\begin{enumerate}
			\item 
			$$\nu_p((\mathbb{Z}_K:\mathbb{Z}[\alpha]))\geq\sum_{i=1}^{r}\text{ind}_{\phi_i}(F).$$  
			The equality holds if $F(x) \text{ is }p$-regular.
			\item 
			If $F(x) \text{ is }p$-regular, then
			$$p\mathbb{Z}_K=\prod_{i=1}^r\prod_{j=1}^{t_i}\prod_{s=1}^{s_{ij}}\mathfrak{p}_{ijs}^{e_{ij}}$$
			is the factorization of $p\mathbb{Z}_K$ into powers of prime ideals of $\mathbb{Z}_K$, where $e_{ij}$ is the smallest positive integer satisfying $e_{ij}\la_{ij}\in \Z$ and the residue degree of $\mathfrak{p}_{ijs}$ over $p$ is given by $f_{ijs}=\deg(\phi_i)\times \deg(\psi_{ijs})$ for every $(i,j,s)$.
		\end{enumerate}
	\end{thm}
	\smallskip
    When  the theorem  of Ore fails, that is $F(x)$ is not $p$-regular, then in order to complete the factorization of $F(x)$, Guardia, Montes, and Nart introduced the notion of {\it high order Newton polygon}.  By analogous to the first order, for each order $r$, the authors  of  \cite{GMN} introduced the valuation $\om_2$ of order $r$, the key polynomial $\ph_2(x)$ of such a valuation,  $N_r(F)$ the Newton polygon of any polynomial $F(x)$ with respect to  $\om_2$ and $\ph_2(x)$, and for every side of $N_r(F)$ the residual polynomial $R_r(F)$, and the index of $F(x)$ in order $r$.  For more details, we refer to \cite{GMN}.
									\section{proofs of our main results}
					{\it Proof of Theorem \ref{intclos}}.\\	
Let  $p$ be a prime integer. Since $\triangle=-(6^6a^7+7^7b^6)$ is the discriminant of $F(x)$, then  a necessary condition for   $p$ to be  a  candidate divisor of   $(\Z_K:\Z[\al])$ is that   $p^2$ divides $6^6a^7+7^7b^6$. 
We conclude that  if  $p=2,3$, then  $p$ divides $b$ is  a necessary condition for   $p$ to be  a  candidate divisor of   $(\Z_K:\Z[\al])$. If $p=7$, then $7$ divides $a$ is  a necessary condition for   $p$ to be  a  candidate divisor of   $(\Z_K:\Z[\al])$. If  $p\not\in \{2,3,7\}$, then  $p^2$ divides $6^6a^7+7^7b^6$ is  a necessary condition for   $p$ to be  a  candidate divisor of   $(\Z_K:\Z[\al])$.  
\begin{enumerate}
\item
If $p$ divides $a$ and $b$, then $p$ does not divide  $(\Z_K:\Z[\al])$ if and only if $\nu_p(b)=1$.
\item
 For $p=2$,
  \begin{enumerate}
 \item
 If $2$ divides $b$ and does not $a$, then  $\ol{F(x)}=x(x-1)^2(x^2+x+1)^2$. Let $\ph_1=x-1$ and $\ph_2=x^2+x+1$. 
 Since $F(x)=\dots-(4x-2)\ph_2+(a+1)x+b$ and $F(x)=\dots+(a+7)\ph_1+(a+1+b)$, we conclude that  $2$ does not divide  $(\Z_K:\Z[\al])$ if and only if $\nu_2((a+1)x+b)=1$ and $\nu_2(b+a+1)=1$, which means $b\equiv 2 \md4$ and $a\equiv 3 \md4$ or $b\equiv 0 \md4$ and $a\equiv 1 \md4$.
  \end{enumerate}
 \item
For $p=3$,
  \begin{enumerate}
 \item
 If $3$ divides $b$ and  $a\equiv 1\md3$, then  $\ol{F(x)}=x(x^2+1)^3$. Let $\ph=x^2+1$. 
 Since $F(x)=x\ph^3-3x\ph^2+3x\ph+(a-1)x+b$, we conclude that  $3$ does not divide  $(\Z_K:\Z[\al])$ if and only if $\nu_3((a-1)x+b)=1$, which means that $a\not\equiv 1 \md9$ or  $b\not\equiv 0 \md9$.
 \item
If $3$ divides $b$ and  $a\equiv -1\md3$, then  $\ol{F(x)}=x(x-1)^3(x+1)^2$. Let $\ph_1=x-1$ and $\ph_2=x+1$. 
 Since  $F(x)=\ph_1^7+7\ph_1^6+21\ph_1^5+35\ph_1^4+35\ph_1^3+21\ph_1^2+(a+7)\ph_1+(a+1+b)$ and $F(x)=\ph_1^7-7\ph_2^6+21\ph_2^5-35\ph_2^4+35\ph_2^3-21\ph_2^2+(a+7)\ph_2+(b-a-1)$, we conclude that  $3$ does not divide  $(\Z_K:\Z[\al])$  if and only if $\nu_2(a+1+b)=1$ and $\nu_3(b-a-1)=1$.
  \end{enumerate}
  \item
  For $p=7$, $7$ divides $a$ and $7$ does not divide $b$,  we have $\ol{F(x)}=(x+b)^7$. Let $\ph=x+b$. 
  $F(x)=\ph^7-7b\ph^6+21b^2\ph^5-35b^3\ph^4+35b^4\ph^3-21b^5\ph^2+(a+7b^6)\ph+(b-ab-b^7)$, we conclude that  $7$ does not divide  $(\Z_K:\Z[\al])$ if and only if $\nu_7(1-a-b^6)=1$.
\item
For $p\not\in\{2,3,7\}$ such that  $p$ does not divide both $a$ and $b$ and   $p^2$ divides   $6^6a^7+7^7b^6$, let  $t$ be an integer  such that  $6at\equiv -7b \md{p^2}$. Then 
  $(6a)^6F'(t)=7(-7b)^6+6^6a^7\equiv 0\md{p^2}$ and $(6a)^7F(t)\equiv 0\md{p^2}$. Thus by Dedekind's criterion $p$ divides the index $(\Z_K:\Z[\al])$.
  \end{enumerate}
	\smallskip
	
	{For the proofs of Theorems \ref{npib2} and \ref{npib3}, we need the following lemma, which  characterizes the prime common index divisors of $K$}.
	 	\begin{lem}\label{index}
		Let $p$ be a rational prime integer and $K$ be a number field. For every positive integer $f$, let $\mathcal{P}_f$ be the number of distinct prime ideals of $\mathbb{Z}_K$ lying above $p$ with residue degree $f$ and $\mathcal{N}_f$ the number of monic irreducible polynomials of $\mathbb{F}_p[x]$ of degree $f$.  Then $p$ is a prime common index divisor of $K$ if and only if $\mathcal{P}_f>\mathcal{N}_f$ for some positive integer $f$.
	\end{lem}
	{\it Proof of Theorem \ref{npib2}}.\\
Since $\triangle=-(6^6a^7+7^7b^6)$ is the discriminant of $F(x)$, if $2$  does not divide $b$, then $2$ does not divide $(\Z_K:\Z[\al])$. Thus $\nu_2(i(K))=0$. Recall that according to Engstrom's results, in order to evaluate $\nu_2(i(K))$, we need to factorize $2\Z_K$ into powers of prime ideals of $\Z_K$ lying above $2$. The following table provides the factorization of $2\Z_K$. Now we assume that  $2$   divides $b$. Then we have the following cases:
	\begin{enumerate}
		\item		
		If $2$ divides  $a$ and $b$, then for $\ph = x$, we have  $ \ol{F(x)}= \ph^7$ in $\F_2[x]$.
		\begin{enumerate}
		\item 
		If $\nph{F}=S$ has a single side, that is  $\nu_2(a)\ge \nu_2(b)$, then  the side $S$ is of degree $1$. Thus  there is a unique prime ideal of $\Z_K$ lying above $2$.
		 \item 
		If $\nph{F}=S_1+S_2$ has two sides joining $(0,\nu_2(b))$, $(1,\nu_2(a))$, and $(7,0)$, that is $\nu_2(a)+1\le \nu_2(b)$, then $S_1$ is of degree $1$, and so it provides a unique prime ideal of  $\Z_K$ lying above $2$ with residue degree $1$. Let $d$ be the the degree of $S_2$.  
		\begin{enumerate}
		\item 
		If $\nu_2(a)\not\in\{ 2,4\}$, then  $S_2$ is of degree $1$, and so there are exactly two prime ideals of $\Z_K$ lying above $2$ with residue degree $1$ each.
		\item
		If $\nu_2(a)=2$,   then   the slope of $S_2$ is $\frac{-1}{3}$ and $R_{1}(F)(y)=(y+1)^2$ is the residual polynomial of $F(x)$ attached to $S_2$. Thus we have to use second order Newton polygon. Let $\om_2=3[\nu_2, x, 1/3]$ be the valuation of second order Newton polygon and  $\ph_2$ the key polynomial of $\om_2$, where $[\nu_2, x, 1/3]$ is the augmented valuation of $\nu_2$ with respect to $x$ and $1/3$. Let $N_2(F)$ the $\ph_2$-Newton polygon of $F(x)$ with respect to the valuation $\om_2$. It follows that:\\  
		If $\nu_2(b)=3$, then for $\ph_2=x^3+2x+2$, we have $F(x)=x\ph_2^2+(4-4x-4x^2)\ph_2+8x^2+(a-4)x+b-8$. It follows that if $\nu_2(a-4)=3$, then $N_2(F)=T$ has a single side joining $(0,10)$ and $(2,7)$. Thus $T$ is of degree $1$ and so $S_2$ provides a unique prime ideal of $\Z_K$ lying above $2$.
		If $\nu_2(a-4)\ge 4$ and $\nu_2(b-8)\ge 4$, then  $N_2(F)=T$ has a single side joining $(0,11)$, $(1,9)$ and $(2,7)$, with $R_2(F)(y)=y^2+y+1$, which is irreducible over $\F_2=\F_0$. Thus $S_2$ provides a unique prime ideal of $\Z_K$ lying above $2$ with residue degree $2$. Hence $2$ is not a common index divisor of $K$.\\
				If $\nu_2(b)\ge 4$ and $\nu_2(a+4)=3$, then  for $\ph_2=x^3+2$,  we have $F(x)=x\ph_2^2-4x\ph_2+(a+4)x+b$ is the $\ph_2$-expansion of $F(x)$, and so $N_2(F)=T$ has a single side joining $(0,10)$ and $(2,7)$. In this case the side $T$ is of degree $1$ and $S_2$ provides a unique prime ideal of $\Z_K$ lying above $2$.
		If $\nu_2(b)= 4$ and $\nu_2(a+4)\ge 4$, then  , then  for $\ph_2=x^3+2$, $N_2(F)=T$ has a single side joining $(0,12)$ and $(2,7)$. Thus $T$ is of degree $1$ and $S_2$ provides a unique prime ideal of $\Z_K$ lying above $2$.\\				If $\nu_2(b)\ge 5$ and $\nu_2(a+4)=4$, then   for $\ph_2=x^3+2$,  we have  $F(x)=x\ph_2^2-4x\ph_2+(a+4)x+b$ be the $\ph_2$-expansion of $F(x)$ and 
 $N_2(F)=T$ has a single side joining $(0,13)$, $(1,10)$ and $(2,7)$. So $T$ is of degree $2$ with attached residual polynomial $R_2(F)=y^2+y+1$ irreducible over $\F_2=\F_0$. Thus $S_2$ provides a unique prime ideal of $\Z_K$ lying above $2$ with residue degree $2$.\\
		If $\nu_2(b)\ge 5$ and $\nu_2(a+4)\ge 5$, then   for $\ph_2=x^3+2$, 
 $N_2(F)=T_1+T_2$ has two sides joining $(0,v)$ , $(1,10)$ and $(2,7)$ with $v\ge 15$. So each $T_i$ has  degree $1$ and $S_2$ provides two prime ideals of $\Z_K$ lying above $2$ with residue degree $1$ each. As $S_1$  provides a  prime ideal of  $\Z_K$ lying above $2$ with residue degree $1$, we conclude that there are three  prime ideals of $\Z_K$ lying above $2$ with residue degree $1$ each, and so $2$ is a common index divisor of $K$. In this last case, $2\Z_K=\p_{111}\p_{121}^3\p_{131}^3$ with residue degree $1$ each factor. Based on Engstrom's result, we conclude that $\nu_2(i(K))=1$.
\item
		Similarly for $\nu_2(a)=4$,  $\nu_2(i(K))\ge 1$
  If and only if $\nu_2(b)\ge 7$ and $\nu_2(a+16)\ge 7$. In this  case, $2\Z_K=\p_{111}\p_{121}^3\p_{131}^3$ with residue degree $1$ each factor. Based on Engstrom's result, we conclude that $\nu_2(i(K))=1$.
		\end{enumerate} 
  \end{enumerate} 
\item
$a\equiv 1\md2$ and $b\equiv 0\md2$. In this case $\ol{F(x)}=x(x-1)^2(x^2+x+1)^2$ modulo $2$.
Let $\ph=x-1$,  $g(x)=x^2+x+1$, $F(x)=\dots-21\ph^2+(7+a)\ph+(b+a+1)$, and $F(x)=(x-3)g^3+(5+3x)g^2-(4x+2)g+(a+1)x+b$.
 Since $x$ provides a unique  prime ideal of $\Z_K$ lying above $2$, we conclude that   $2$ is a common index divisor of $K$ if and only if $\ph$ provides  two  prime ideals of $\Z_K$    
lying above $2$ of degree $1$ each or $\ph$ provides  a unique  prime ideal of $\Z_K$    
lying above $2$ of degree $2$ and $g$ provides  at least one   prime ideal of $\Z_K$    
lying above $2$ of degree $2$ or also $g$ provides  two  prime ideals of $\Z_K$    
lying above $2$ of degree $2$ each. That is if and only if  one of the following conditions holds:
\begin{enumerate}
 \item
If $a\equiv 1\md4$ and $b\equiv 2\md4$, then $\nu_2(\triangle)\ge 7$  and  $N_g^+(F)$ has a single  side of height $1$, and so $g$ provides a unique prime ideal $\p_{311}$ of $\Z_K$ lying above $2$ with residue degree $2$. For $\npp{F}$, let $u=\frac{-7b_2}{3a}$. Then $u\in\Z_2$. Let $F(x+u)=x^7+\dots+21u^5x^2+Ax+B$, where $A=7u^6+a=\frac{-\triangle}{6^6a^6}$ and $B=u^7+au+b=\frac{b\triangle}{6^7a^7}$. 
It follows that $\nu_2(A)=\nu_2(B)=\nu_2(\triangle)-6$, and so  $\npp{F}=S_1$ has a single side joining $(0, \nu_2(\triangle)-6)$ and $(2,0)$. Thus, if  $\nu_2(\triangle)$ is odd, then $\ph$ provides  a unique prime ideal $\p_{211}$ of $\Z_K$ lying above $2$ with residue degree $1$.
If  $\nu_2(\triangle)=2(k+3)$  for some positive integer  $k$, then let $F(x+u+2^k)=x^7+\dots+21(u+2^k)^5x^2+A_1x+B_1$,
where $A_1=7u^6+a+3\cdot 2^{k+1}u^5+2^{2k}D=A+3\cdot 2^{k+1}u^5+2^{2k}D$  and $B_1=B+A\cdot 2^k+2^{2k}\cdot 21u^5+ 2^{3k}H=2^{2k}(\frac{b_2\triangle_2}{3^7a^7}+21u^5)+ 2^{3k}H$ for some  $D\in \Z_2$ and $H\in \Z_2$. Thus,  $B_1=2^{2k}(3\cdot a\cdot b_2\triangle_2+3\cdot a\cdot b_2)+ 2^{2k+3}L$  for some $L\in \Z_2$. Hence if $k\ge 2$, then  $\nu_2(A_1)=k+1$ and $\nu_2(B_1)\ge 2k+1$. More precisely,  if $\triangle_2\equiv 1 \md4$, then $\nu_2(B_1)=2k+1$, and so $\ph$ provides a  unique prime ideal $\p_{211}$ of $\Z_K$ lying above $2$ with residue degree $1$.   If $\triangle_2\equiv 3 \md4$, then $\nu_2(B_1)\ge 2k+2$. It follows that if  $\nu_2(B_1)= 2k+2$, and so $\ph$ provides a  unique prime ideal $\p_{211}$ of $\Z_K$ lying above $2$ with residue degree $2$. If  $\nu_2(B_1)\ge 2k+3$, then $\nu_2(B_1)\ge 2k+3$, and so $\ph$ provides two  prime ideals  of $\Z_K$ lying above $2$ with residue degree $1$ each. In these  last two cases, we have $2$ divides $i(K)$ and $\nu_2(i(K))=1$.\\
For $k=1$, we have  $\nu_2(\triangle)=8$ and  $a\equiv 5\md8$. 
 In this case $\ol{F(x)}=x(x-1)^2(x^2+x+1)^2$ modulo $2$.
Let $\ph=x-1$,  $g(x)=x^2+x+1$, $F(x)=\dots-21\ph^2+(7+a)\ph+(b+a+1)$, and $F(x)=(x-3)g^3+(5+3x)g^2-(4x+2)g+(a+1)x+b$. Since $x$ provides a unique prime ideal of of $\Z_K$ lying above $2$ with residue degree $1$ and $g$ provides a unique prime ideal of of $\Z_K$ lying above $2$ with residue degree $2$,  we conclude that  $\nu_2(i(K))\ge 1$ if and only if   $\ph$ provides a unique prime ideal of of $\Z_K$ lying above $2$ with residue degree $2$ or  $\ph$ provides two distinct prime ideals of of $\Z_K$ lying above $2$ with residue degree $1$ each.   If  $(a,b)\in \{(5,10), (13,2)\}\md{16}$, then $\ph$ provides a unique prime ideal of $\Z_K$ lying above $2$ with residue degree $2$ and so $\nu_2(i(K))=1$. If  $(a,b)\in\{(5,10), (13,10)\}\md{16}$, then $\ph$ provides a unique prime ideal of $\Z_K$ lying above $2$ with residue degree $1$ and so $\nu_2(i(K))=0$. For   $(a,b)\in\{(5,6), (5,14),(13,6), (13,14)\}\md{16}$, let us replace $\ph$ by $\ph'=x-3$ and consider the $\ph'$-Newton polygon of $F(x)$ with respect to $\nu_2$. It follows that 
If  $(a,b)\in\{(5,6), (13,14)\}\md{16}$, then $\ph'$ provides two prime ideals of $\Z_K$ lying above $2$ with residue degree $1$ each and so $\nu_2(i(K))=1$. If  $(a,b)\in\{(5,14), (13,6)\}\md{16}$, then $\ph'$ provides a unique prime ideal of $\Z_K$ lying above $2$ with residue degree $1$ and so $\nu_2(i(K))=0$.
\item
$a\equiv 3\md4$ and    $b\equiv -(a+1)\md8$ because $\npp{F}$ has two sides.
 \item
  If $a\equiv 3\md8$ and   $b\equiv 0\md8$, then $\ph$ provides a unique  prime ideal of $\Z_K$ lying above $2$ with residue degree $2$ and $g$ provides two  prime ideals of $\Z_K$ lying above $2$ with residue degree $2$ each because $N_g^+(F)$ has  a single side of degree $2$ with $(1+x)y^2+y+x=(x+1)(y-1)(y-x)$ its attached residual polynomial of $F(x)$.  In this case $2\Z_K=\p_{111}\p_{211}\p_{221}\p_{311}\p_{312}$ with residue degrees $f_{111}=1$ and $f_{211}=f_{311}=f_{312}=2$, and so $\nu_2(i(K))=3$.
\item
 $a\equiv 7\md8$ and    $b\equiv 0\md8$. In this case $\ph$ provides a unique  prime ideal of $\Z_K$ lying above $2$ with residue degree $2$  and $N_g^+(F)$ has two sides.
  More precisely, $2\Z_K=\p_{111}\p_{211}\p_{221}\p_{311}\p_{312}$ with residue degrees $f_{111}=f_{211}=f_{221}=1$ and $f_{211}=f_{311}=f_{312}=2$, and so $\nu_2(i(K))=3$.
 \item
If $a\equiv 5\md8$ and    $b\equiv -(a+1)\md{16}$ because if $b\equiv -(a+1)\md{32}$, then $\npp{F}$ has two sides   and if $b\equiv -(a+1)+16\md{32}$, then  $\npp{F}$ has a single side of degree $2$, which provides a single prime ideal of $\Z_K$ lying above $2$ with residue degree $2$ and $N_g^+(F)$ has  a single side of degree $1$. Thus there are $2$ prime ideals of $\Z_K$ lying above $2$ with residue degree $2$ each. 
\item
If $a\equiv 5\md8$ and    $\nu_2(b+(a+1))=2$, then $\nu_2(b-(a+1))\ge 3$. If $\nu_2(b-(a+1))= 3$, then for  $\ph=x+1$, we have $\npp{F}$ has a single side of degree $1$.
Since $\nu_2(a+1)=1$, then $N^+_g(F)$ has a single side of height $1$. Thus there are two prime ideals of $\Z_K$ lying above $2$ with residue degree $1$ each and one prime ideal with residue degree $2$.  If $\nu_2(b-(a+1))= 4$, then for  $\ph=x+1$, we have $\npp{F}$ has a single side of degree $2$ and its attached residual polynomial of $F$ is $R_1(F)(y)=y^2+y+1$. Since $b\equiv 6\md{8}$, we conclude that $N_g^+(F)$ has  a single side of degree $1$, then there are $2$ prime ideals of $\Z_K$ lying above $2$ with residue degree $2$ each, and so $2$ divides $i(K)$. If $\nu_2(b-(a+1))\ge 5$, then for  $\ph=x+1$, we have $\npp{F}$ has two  sides of degree $1$ each, and so there are $3$ prime ideals of $\Z_K$ lying above $2$ with residue degree $1$ each, and so $2$ divides $i(K)$.
		\end{enumerate} 
\end{enumerate}
				
{\it Proof of Theorem \ref{npib3}}.\\
Since $\triangle=-(6^6a^7+7^7b^6)$ is the discriminant of $F(x)$, by the formula {$\triangle=(\Z_K:\Z[\al])^2d_K$},   if $3$ does not divide $b$, then $3$ does not divide $(\Z_K:\Z[\al])$, and so $\nu_3(i(K))=0$. Now assume that $3$  divides $b$. Recall that according to Lemma \ref{index}, $\nu_3(i(K))\ge 1$ if and only if there are at least four primes ideal of $\Z_K$  lying above $3$ with residue degree $1$ each.
We have the following cases:
\begin{enumerate}
		\item 
		If  $3$   divide $a$ and $b$, then for  $\ph=x$, $\ol{F(x)}= \ph^7$ in $\F_3[x]$. It follows that:
		\begin{enumerate}
		\item 
		If $\nu_3(a)\ge \nu_3(b)$, then $\nph{F}$ has a single side of degree $1$, and so there is a unique prime ideal of $\Z_K$ lying above $3$.
		 \item 
		If $\nu_3(a)+1\le \nu_3(b)$, then $\nph{F}=S_1+S_2$ has two sides joining $(0,\nu_3(b))$, $(1,\nu_3(a))$, and $(7,0)$. Since $S_1$ is of degree $1$, $S_1$ provides exactly one prime ideal of $\Z_K$  lying above $3$ with residue degree $1$. Thus  $\nu_3(i(K))\ge 1$ if and only if $S_2$ provides  at least three primes ideal of $\Z_K$  lying above $3$ with residue degree $1$ each.
    If $\nu_3(a)\not\in\{2,3, 4\}$, then  $S_2$ is of degree $1$, and so $S_2$ provides exactly one prime ideal of $\Z_K$ lying above $3$ with residue degree $1$ each.
		If $\nu_3(a)\in\{2,4\}$,  then $S_2$ is of degree $2$, and so $S_2$ provides at most two  prime ideal of $\Z_K$ lying above $3$. Hence $3$ is not  a common index divisor of $K$.
If  $\nu_3(a)=3$,  then $S_2$ is of degree $3$ and its attached residual polynomial of $F(x)$ is $R_1(F)(y)=y^3+a_3=(y+a_3)^3$. So we have to use second order Newton polygon. Let $\om_2=2[\nu_3, x, 1/2]$ be the valuation of second order Newton polygon, $\ph_2=x^2+3a_3$ the key polynomial of $\om_2$ and $N_2(F)$ the $\ph_2$-Newton polygon of $F(x)$ with respect to  $\om_2$. It follows that:
If $a_3\equiv 1\md3$, then  for $\ph_2=x^2+3$, we have $F(x)=x\ph_2^3-9x\ph_2^2+27x\ph_2+(a-27)x+b$ is the $\ph_2$-expansion of $F(x)$. We have the following cases:
\begin{enumerate}
\item
		If $\nu_3(b)=4$, then $N_2(F)=T$ has a single side joining $(0,8)$ and $(3,7)$. Thus $T$ is of degree $1$ and $S_2$ provides a unique prime ideal of $\Z_K$ lying above $3$ with residue degree $1$. 
\item
		If $\nu_3(b)\ge 5$ and $\nu_3(a-27)=4$, then $N_2(F)=T$ has a single side joining $(0,9)$ and $(3,7)$. Thus $T$ is of degree $1$ and $S_2$ provides a unique prime ideal of $\Z_K$ lying above $3$ with residue degree $1$.
\item
		If $\nu_3(b)= 5$ and $\nu_3(a-27)\ge 5$, then $N_2(F)=T$ has a single side joining $(0,10)$ and $(3,7)$ and its  attached residual polynomial of $F$ is $R_2(F)(y)=xy^3+xy+b_3$, which is irreducible over $\F_2=\fph$ because $\ph$ is of degree $1$. Thus  $S_2$ provides a unique prime ideal of $\Z_K$ lying above $3$ with residue degree $3$.  
		\item
		If $\nu_3(b)\ge 6$ and $\nu_3(a-27)\ge 5$, then $N_2(F)=T_1+T_2$ has two  sides joining $(0,v)$, $(2,9)$ and $(3,7)$ with $v\ge 11$. Thus $T_1$ is of degree $1$, $T_2$ of degree $2$  and $R_2(F)(y)=xy^2+x$ is its attached residual polynomial of $F(x)$, which is irreducible over $\F_2=\fph$. Thus  $S_2$ provides a unique prime ideal of $\Z_K$ lying above $3$ with residue degree $1$ and a unique prime ideal of $\Z_K$ lying above $3$ with residue degree $2$.
		\end{enumerate}
  Similarly, for $a_3\equiv -1\md3$, let $\ph_2=x^2-3$. Then  $F(x)=x\ph_2^3+9x\ph_2^2+27x\ph_2+(a+27)x+b$ is the $\ph_2$-expansion of $F(x)$. By analogous to the case $a_3\equiv 1\md3$, in every case $3$ does not divide $i(K)$.
		If $a\equiv 1\md3$, then $\ol{F(x)}=x(x^2+1)^3$ in $\F_3[x]$. Then there are exactly a unique  prime ideal of $\Z_K$ lying above $3$ with residue degree $1$ and the others prime ideals of $\Z_K$ lying above $3$ are of residue degrees at least $2$ each. Hence $\nu_3(i(K))=0$.
  \item
		If $a\equiv -1\md3$, then $\ol{F(x)}=x(x-1)^3(x+1)^3$ in $\F_3[x]$.
		 		  Let $\ph_1=x-1$, $\ph_2=x+1$, 				
				 $F(x)=\ph_1^7+7\ph_1^6+21\ph_1^5+35\ph_1^4+35\ph_1^3+21\ph_1^2+(7+a)\ph_1+(b+a+1)$, and 
								$F(x)=\ph_2^7-7\ph_2^6+21\ph_2^5-35\ph_2^4+35\ph_2^3-21\ph_2^2+(7+a)\ph_2+(b-(a+1))$.  
It follows that:
\begin{enumerate}
\item
If  $b\equiv -(1+a) \pm 3 \md{9}$ and  $b\equiv (1+a) \pm 3 \md{9}$, then for $N_{\ph_i}^+(F)$ has a single side of height $1$ for every $i=1,2$. Thus there are exactly three prime ideals of $\Z_K$ lying above $3$ with residue degree $1$ each. 
\item
If $a\equiv 8\md{9}$ and  $b\equiv 0 \md{9}$, then $\nu_3(b+(1+a))\ge 2$ and $\nu_3(b-(1+a))\ge 2$. Thus  $x$ a provides a unique prime ideal of $\Z_K$ lying above $3$ with residue degree $1$,  and each of  $\ph_1$ and $\ph_2$ provides two  prime ideals of $\Z_K$ lying above $3$ with residue degree $1$ each. 
 In this two cases $\nu_3(i(K))=2$.
\item
If $a\equiv 5\md{9}$ and  $b\equiv 3 \md{9}$, then  $\nu_3(b+(1+a))\ge 2$ and $\nu_3(b-(1+a))=1$. Thus  each of $x$ and  $\ph_2$ provides a unique prime ideal of $\Z_K$ lying above $3$ with residue degree $1$,  and   $\ph_1$ provides two  prime ideals of $\Z_K$ lying above $3$ with residue degree $1$ each. 
Similarly, if $a\equiv 5\md{9}$ and  $b\equiv 6 \md{9}$, then  $\nu_3(b-(1+a))\ge 2$ and $\nu_3(b+(1+a))=1$. Thus each of $x$ and  $\ph_1$ provides a unique prime ideal of $\Z_K$ lying above $3$ with residue degree $1$,   $\ph_2$ provides two  prime ideals of $\Z_K$ lying above $3$ with residue degree $1$ each.  In these two cases $\nu_3(i(K))=1$.
\item
If 	$a\equiv 2\md9$ and $\nu_3(b)\ge 2$, then $\nu_3(b-(a+1))=1$ and $\nu_3(b+(a+1))=1$. Thus $3$ does not divide $(\Z_K:\Z[\al])$.
\item
If 	$a\equiv 2\md9$ and $b\equiv (1+a) \pm 9 \md{27}$, then $N_{\ph_2}^+(F)$ has a single  side joining $(0,2)$ and $(3,0)$ and $N_{\ph_1}^+(F)$ has  a single side joining $(0,1)$ and $(3,0)$. Thus there are $3$ prime ideals of $\Z_K$ lying above $3$ with residue degree $1$ each, and so  $\nu_3(i(K))=0$.
\item
Similarly, if 	$a\equiv 2\md9$ and $b\equiv -(1+a) \pm 9 \md{27}$, then  there are $3$ prime ideals of $\Z_K$ lying above $3$ with residue degree $1$ each, and so  $\nu_3(i(K))=0$.
\item
If 	$a\equiv 2\md9$ and $\nu_3(b)=1$, then $\nu_3(\triangle)\ge 8$. Let 
 $u=\frac{-7b_3}{2a}$. Then $u\in\Z_3$. Let $\ph=x-u$ and $F(x+u)=x^7+\dots+35u^4x^3+21u^5x^2+Ax+B$, where $A=7u^6+a=\frac{-\triangle}{6^6a^6}$ and $B=u^7+au+b=\frac{b\triangle}{6^7a^7}$. It follows that $\nu_3(A)=\nu_3(B)=\nu_3(\triangle)-6$, and so  $\npp{F}=S_1$ has a single side joining $(0, \nu_3(\triangle)-6)$ and $(2,0)$. Remark that since  $\nu_3(b)=1$ and $\nu_3(B)\ge 2$,  $\nu_3(-u^7-au+b)=1$, and so $(x+u)$ provides a unique prime  ideal of $\Z_K$ lying above $3$ with residue degree $1$. Thus  $\nu_3(i(K))\ge 1$  if and only if $\ph$ provides at least two
 prime ideals of $\Z_K$ lying above $3$ with residue degree $1$ each.
  Let $v=\nu_3(\triangle)$.  We have the following cases:
  \begin{enumerate}
\item
If $v=8$, then $\npp{F}$ has a single side of degree one, and so $\ph$ provides a unique prime ideal of $\Z_K$ lying above $3$ with residue degree $1$.
\item
If $v=9$, then  $\npp{F}=S$ has a single side joining $(0,3)$ and $(3,0)$ with $R_1(F)(y)=-u^4y^3+u^5y^2+B_3$ its attached residual polynomial of $F(x)$. Since $a\equiv -1\md3$ and $B=\frac{b\triangle}{6^7a^7}$, we have $u\equiv -b_3\md3$ and $B_3\equiv b_3\triangle_3\md3$. Thus $R_1(F)(y)=-y^3-b_3y^2+ b_3\triangle_3$.  Since $R_1(F)(y)$ is square free and $R_1(F)(0)\neq 0$, then $R_1(F)(y)$ has at most one root in $\fph$. Thus $S$ provides at most  a unique prime ideal of $\Z_K$ lying above $3$ with residue degree $1$. Therefore, $\nu_3(i(K))=0$.
\item
If $v\ge 10$, then  $\npp{F}=S_1+S_2$ has two sides joining $(0,v-6)$, $(1,v-6)$ and $(3,1)$.
It follows that
Since  $S_2$ is of degree $1$,  it provides a unique prime ideal of $\Z_K$ lying above $3$ with residue degree $1$. Moreover,
if $v$ is even then $S_1$ is of degree $1$, and so $\ph$ provides two prime ideals of $\Z_K$ lying above $3$ with residue degree $1$ each. In this case $\nu_3(i(K))=1$.  
  If $v=2(k+3)+1$, then $S_1$ is of degree $2$ with residual polynomial $R_1(F)(y)=uy^2+b_3\triangle_3$. Since $a\equiv -1\md3$, we have $2a\equiv 1\md3$ and $u\equiv -b_3\md3$. Thus $R_1(F)(y)=-b_3(y^2-\triangle_3)$. It follows that if $(\frac{\triangle_3}{3})=1$, then $R_1(F)(y)$ has two different factors of degree $1$ each, and so $S_1$ provides two  prime ideals of $\Z_K$ lying above $3$ with residue degree $1$ each. In this case there are exactly five  prime ideals of $\Z_K$ lying above $3$ with residue degree $1$ each and according to Engstrom's results $\nu_3(i(K))=2$. But if $(\frac{\triangle_3}{3})=-1$, then $R_1(F)(y)$  is irreducible over $\fph=\F_3$, and so $S_1$ provides a unique  prime ideal of $\Z_K$ lying above $3$ with residue degree $2$. In this last case there are exactly three prime ideals of $\Z_K$ lying above $3$ with residue degree $1$ each, and  so $\nu_3(i(K))=0$.
 				\end{enumerate}
     \end{enumerate}
 				\end{enumerate}
     \end{enumerate}
 {\it Proof of Theorem \ref{pge5}}.\\
	
	    We start by showing that $5$ does not divide $i(K)$ for every integers of $a$ and $b$ such that $x^7+ax+b$ is irreducible.
	    \begin{enumerate}
	     \item 
	    If  $5^2$ does not divide $\triangle$, then $5$ does not divide $(\Z_K:\Z[\al])$. Consequently, $\nu_5(i(K))=0$.
	    \item
	    Now assume that   $5^2$  divides $\triangle=-(6^6a^7+7^7b^6)$. Then $6^6a^7+7^7b^6\equiv 0\md5$. Since $a^5\equiv a\md5$ and $b^5\equiv b\md5$, then $a^2\equiv 3b\md5$, which means $(a,b)\in\{(0,0), (2,2), (3,2), (1,3), (4,3)\}\md5$. In order to show that $\nu_5(i(K))=0$ it suffices to show that  for every value $(a,b)\in\Z^2$ such that $x^7+ax+b$ is irreducible and $(a,b)\in \{(0,0), (1,2), (-1,2), (2,3), (3,3)\}\md5$   there are at most four prime ideals of $\Z_K$ lying above $5$ with residue degree $1$, where $K$ is the number field generated by a complex root of $x^7+ax+b$ . 
	    \begin{enumerate}
	    \item 
	    For $(a,b)\equiv (0,0)\md5$, if $\nu_5(a)\ge \nu_5(b)$, then  $\nph{F}=S$ has a single side and it is  of degree $1$. Thus there is a unique prime ideal $\p$ of $\Z_K$ lying above $5$ with residue degree $1$. More precisely $5\Z_K=\p^7$. \\
     If  $\nu_5(a)+1\le \nu_5(b)$, then $\nph{F}=S_1+S_2$ has two sides. More precisely, $S_1$ is of degree $1$. Let $d$ be degree of $S_2$. Since $6$ is the length of $S_2$, then $d\in\{1,2,3\}$. Thus $S_1$ provides a unique prime ideal $\p$ of $\Z_K$ lying above $5$ with residue degree $1$ and $S_2$ provides at most three prime ideals $\p$ of $\Z_K$ lying above $5$ with residue degree $1$ each.
     	    \item
	    For $(a,b)\equiv (2,2)\md5$, since  $\ol{F(x)}=(x-1)(x-3)^2(x^4+2x^3+4x^2+2x+2)$ in $\F_5[x]$, there is at most three  prime ideals $\p$ of $\Z_K$ lying above $5$ with residue degree $1$.
      \item
	    For $(a,b)\in\{ (3,2),  (4,3)\}\md5$, since  $\ol{F(x)}$ is irreducible over  $\F_5$, there is a unique  prime ideal $\p$ of $\Z_K$ lying above $5$ with residue degree $7$.
      \item
	    For $(a,b)\equiv (1,3)\md5$, since  $\ol{F(x)}=(x-1)(x^6+x^5+x^4+x^3+x^2+x+2)$ in $\F_5[x]$, there is exactly a unique  prime ideals $\p$ of $\Z_K$ lying above $5$ with residue degree $1$.
	\end{enumerate}
	   \end{enumerate}
     We conclude that in all cases $\nu_5(i(K))=0$.\\
	For $p\ge 7$, since the field $K$ is of degree $7$, there are at most $7$ prime ideals of $\Z_K$ lying above $p$. The fact that  there at least $p\ge 7$ monic irreducible polynomial of degree $f$ in $\F_p[x]$ for every positive integer $f\in\{1,2,3\}$, we conclude that  $p$ does not divide $i(K)$.
	
	{\it Proof of Proposition \ref{mono}}.\\
		First according to Theorem \ref{intclos} and the hypotheses of Example \ref{mono}, $2$ is the unique prime integer candidate to divide $ind(\al)$.
Let $\ph=x$. Then $\ol{F(x)}=\ph^7$ in $\F_2[x]$ and $\nph{F}=S$ has a single side of degree GCD$(7,v)=1$. Thus  $F(x)$ is irreducible over $\Q_2$. Let $K$ be the number field generated by a  root $\al$ of $F(x)$. Since $F(x)$ is irreducible over $\Q_2$, 
there is a unique valuation $\om$ of $K$ extending $\nu_2$. 
By Theorem \ref{ore}, we have $\nu_2(\Z_K:\Z[\al])\ge ind_\ph(F)\ge 1$, and so $F(x)$ is not a  monogenic polynomial. Let $\theta =\frac{\al^x}{2^y}$, where $(x,y)$ is the unique   solution of integers of the Diophantine equation $vx-7y=1$ and $0\le x \le 6$.
	Then $\theta \in K$. Since $v$ and $7$ are coprime, we conclude that $K=\Q(\theta)$.	 Let us show that $\Z_K=\Z[\theta]$, and so $K$ is monogenic.  By \cite[Corollary 3.1.4]{En}, in order to show that $\theta\in \Z_K$, we need to show that $\om(\theta)\ge 0$, where $\om$ is the unique valuation of $K$ extending $\nu_2$. Since $\nph{F}=S$ has  a single side of slope $-v/7$, we conclude that $\om(\al)=v/7$, and so $\om(\theta)=\frac{xv}{7}-y=\frac{1}{7}$. Let $g(x)$ be the minimal polynomial of $\theta$ over $\Q$. By the formula relating roots and coefficients of a monic polynomial, we conclude that $g(x)=x^{7}+\sum_{i=1}^{7}(-1)^is_ix^{7-i}$, where $s_i=\displaystyle\sum_{k_1<\dots<k_i}\theta_{k_1}\cdots\theta_{k_i}$ and $\theta_{1},\dots, \theta_{7}$ are the $\Q_p$-conjugates of $\theta$. Since there is a unique valuation extending $\nu_2$ to any algebraic extension of $\Q_2$, we conclude that $\om(\theta_i)=1/7$ for every $i=1,\dots,7$. Thus $\nu_2(s_{7})=\om(\theta_{1}\cdots\theta_{7})=7\times 1/7=1$ and $\nu_2(s_{i})\ge i/7$ for every $i=1,\dots, 6$, which  means that $g(x)$ is a $2$-Eisenstein polynomial. Hence $2$ does not divide the index $(\Z_K:\Z[\theta])$. As $2$ is the unique positive prime integer candidate to divide  $(\Z[\al]:\Z[\theta])$, we conclude that for every prime integer $p$, $p$ does not divide $(\Z_K:\Z[\theta])$, which means that $\Z_K=\Z[\theta]$.
\section{Examples} 
Let $F=x^7+ax+b \in \Z[x]$ be a	 monic irreducible polynomial  and $K$ a number field generated by a  root $\al$ of $F(x)$. In the following examples, we calculate the index of the field $K$. First based on Theorem \ref{pge5}, $\nu_p(i(K))=0$ for every prime integer $p\ge 5$. Thus we need only to calculate $\nu_p(i(K))$ for $p=2,3$.
\begin{enumerate}
\item  
For $a=6$ and $b=6$, since ${F(x)}$ is $p$-Eisenstein for every $p=2,3$, we conclude that   $F(x)$ is irreducible over $\Q$, $2$ (resp. $3$) does not divide $(\Z_K:\Z[\al])$. Thus  $2$ (resp. $3$) does not divide $i(K))$, and so   $i(K)=1$.
\item  
For $a=28$ and $b=32$, since $\ol{F(x)}$ is irreducible over $\F_5$,  $F(x)$ is irreducible over $\Q$. By the first item of Theorem \ref{npib2}, we have  $\nu_2(i(K))=1$. By Theorem \ref{npib3},   $\nu_3(i(K))=0$. Thus $i(K)=2$.
\item
For $a=3$ and $b=8$,  $\ol{F(x)}$ is irreducible over $\F_5$,  $F(x)$ is irreducible over $\Q$. Again since $a\equiv 3\md4$ and $b\equiv 0  \md{8}$, by Theorem \ref{npib2}, $\nu_2(i(K))=3$.  By Theorem \ref{npib3},   $\nu_3(i(K))=0$. Thus $i(K)=8$.
\item  
For $a=-1$ and $b=9$,  since $\ol{F(x)}$ is irreducible over $\F_2$,  $F(x)$ is irreducible over $\Q$. Since $2\Z_K$ is a prime ideal of $\Z_K$, $\nu_2(i(K))=0$. Also since $a\equiv 8\md{9}$ and $b\equiv 0\md9$, by Theorem \ref{npib3}, $\nu_3(i(K))=2$. Thus $i(K)=9$.
\item  
For $a=803$ and $b=2112$,   since $\ol{F(x)}$ is irreducible over $\F_5$,  $F(x)$ is irreducible over $\Q$. Since $a\equiv 3\md4$ and $b\equiv 0\md8$, by Theorem \ref{npib2}, $\nu_2(i(K))=3$. Similarly since $a\equiv 5\md9$ and $b\equiv 6\md9$, by Theorem \ref{npib3},  $\nu_3(i(K))=1$. Thus  $i(K)=24$.
\item  
For $a=35$ and $b=72$,   since $\ol{F(x)}$ is irreducible over $\F_{11}$,  $F(x)$ is irreducible over $\Q$. Since $a\equiv 3\md4$ and $b\equiv 0\md8$, by Theorem \ref{npib2}, $\nu_2(i(K))=3$. Similarly since $a\equiv 8\md9$ and $b\equiv 0\md9$, by Theorem \ref{npib3},  $\nu_3(i(K))=2$. Thus  $i(K)=72$.
\end{enumerate}
{\bf Data Availability Statement}\\
Data sharing not applicable to this article as no datasets were generated or analysed during the current study.

\end{document}